\documentclass[10pt,eqno]{article}
\usepackage{latexsym}
\usepackage{amsfonts,amsmath,amssymb}                                                                                                                                                                                                                                                                                                                                                                                                                                                                                                                                                                                      
\usepackage{amsthm}
\usepackage{}
\newtheorem{theorem}{Theorem}[section]
\newtheorem{lemma}{Lemma}[section]
\newtheorem{proposition}{Proposition}[section]

\newtheorem{remark}{Remark}[section]
\usepackage[pdftex]{graphicx}
\usepackage[latin1]{inputenc}
\usepackage{epstopdf}     
\usepackage{setspace}
\usepackage{tikz}
\usepackage[toc]{appendix}
\usepackage{hyperref}

\numberwithin{equation}{section}
\numberwithin{theorem}{section}

\numberwithin{proposition}{section}
\numberwithin{lemma}{section}
\numberwithin{remark}{section}
\numberwithin{equation}{section}    
\begin{document} 

\setcounter{secnumdepth}{3}
\newcommand{\bbox}{\vrule height.6em width.6em
depth0em} 
\parindent0in     
\def\lhead{vv}
\thispagestyle{empty}
\setcounter{page}{1}
\noindent
\def\lhead{vv }
\def\rhead{ Holder}
\thispagestyle{empty}
\setcounter{page}{1}
\noindent
\begin{center}
{\bf {\Large H\"older estimates  of weak solutions to   chemotaxis systems of fast diffusion type.}}
\end{center}

\vspace{1cm}
{\bf 
{ M.Marras \footnote{ Dipartimento di Matematica e Informatica, via Ospedale 72, \!\!09124 Cagliari (Italy),\\
 mmarras@unica.it},
F.Ragnedda \footnote{ Facolt\`a  di Ingegneria e Architettura,  via Marengo 2, 09123 Cagliari, (Italy), ragneddafra11@gmail.com}, 
S.Vernier-Piro \footnote{ Facolt\`a di Ingegneria e Architettura, Via Marengo   2, 09123 Cagliari (Italy), svernier@unica.it}
and V.Vespri\footnote{Dipartimento di Matematica ed Informatica 'U. Dini', viale Morgagni 67/a, 50134 Firenze (Italy)
vincenzo.vespri@unifi.it}  }}\\

\vspace*{1cm}
\thispagestyle{empty}
\setcounter{page}{1}
\noindent

\setcounter{page}{1}
\noindent  
\begin{abstract}
We study a quasilinear chemotaxis system of singular type, where the diffusion operator is given by $\Delta u^m$ with $0<m<1$, corresponding to the fast diffusion regime, and where the chemotactic drift is nonlinear. Since H\"older continuity constitutes the optimal regularity class for weak solutions to the porous medium equation, we establish analogous regularity results for bounded solutions of parabolic--parabolic chemotaxis systems in this setting. The proof is based on a refined De Giorgi--Di Benedetto iteration scheme adapted to the coupled structure of the system. These results advance the understanding of the fine regularity properties of chemotaxis models with nonlinear diffusion, and demonstrate that the interplay between singular diffusion and aggregation exhibits a regularizing mechanism consistent with the porous medium paradigm.
\end{abstract}


\section{Introduction}

We consider the following class of singular chemotaxis systems:
\begin{equation}\label{1.1}
\begin{cases}
u_t = \operatorname{div}(\nabla u^m) - 
\operatorname{div}(
\chi(u,v)u  \nabla v), & \text{in } \mathbb{R}^N \times (t>0), \\[6pt]
v_t = \Delta v - \alpha v + u, & \text{in } \mathbb{R}^N \times (t>0), \\[6pt]
u(x,0) = u_0(x) > 0,\quad v(x,0) = v_0(x) \geq 0, & \text{in } \mathbb{R}^N,
\end{cases}
\end{equation}
with $N \geq 1$, $\frac {(N-2)_+}{N+2}<m<1,$ $\alpha>0,\chi $ a positive function.  The initial data $(u_0,v_0)$ are assumed to satisfy
\begin{equation*}
\begin{cases}
u_0 > 0,\quad u_0 \in L^\infty(\mathbb{R}^N)\cap L^1(\mathbb{R}^N), \quad u_0^m \in H^1(\mathbb{R}^N), \\[6pt]
v_0 \geq 0,\quad v_0 \in L^1(\mathbb{R}^N)\cap W^{1,p}(\mathbb{R}^N),
\end{cases}
\end{equation*}
where $p = \infty$ if $N=1$ and $p= \frac{N(m+1)}{N-m-1}$ if $N\geq 2.$

The system \eqref{1.1} is termed \emph{singular} because the diffusion operator
\begin{align*}
\operatorname{div}(\nabla u^m) = m\, \operatorname{div}(u^{m-1}\nabla u)
\end{align*}
becomes singular at $u=0$ when $0<m<1.$
Here, $u = u(x,t)$ denotes the cell density, while $v = v(x,t)$ represents the concentration of the chemoattractant. The nonlinear diffusion coefficient $D(u)$ and the chemotactic sensitivity $\chi(u,v)$ govern the interplay between dispersive and aggregative effects. In this work, we assume that the chemotactic flux satisfies $\chi(u,v)u =\chi  u^{{\mathfrak q}-1}$ 
with $\chi$ constant  and 
\begin{align}\label{mathfrak q}
1< {\mathfrak q}< \frac{(m+1)(N+2) }{2N },
\end{align}
although any function exhibiting an analogous growth behavior could be considered within the same analytical framework. 
From the modeling viewpoint, it is also relevant to consider the \emph{parabolic-elliptic} reduction, where the second equation of the system is stationary rather than time-evolving  (see \cite{DM}). However, in the present paper we restrict our attention to \emph{parabolic-parabolic} systems, which retain the full time dynamics of both variables and present additional analytical difficulties related to the coupling of the two evolution equations.

\medskip

A pair $(u,v)$ of nonnegative measurable functions defined in $\mathbb{R}^N \times (0,T]$, $T>0$, is called a \emph{local weak solution} of \eqref{1.1} if
\[
u \in C_{\mathrm{loc}}(0,T;L^{m+1}_{\mathrm{loc}}(\mathbb{R}^N)), 
\  |u|^m \in L^2_{\mathrm{loc}}(0,T;W^{1,2}_{\mathrm{loc}}(\mathbb{R}^N)), 
\ v \in L^\infty(0,T;H^1(\mathbb{R}^N)),
\]

and if \eqref{1.1} holds in the weak sense: for every compact ${\cal{K}} \subset \mathbb{R}^N$, every interval $[t_1,t_2]\subset [0,T]$, and every nonnegative test function 
\[
\psi \in W^{1,2}_{\mathrm{loc}}(0,T;L^2({\cal{K}})) \cap L^{2}_{\mathrm{loc}}(0,T; W^{1,2}_{0}
({\cal{K}}))\]
one has
\begin{equation}\label{1.3}
\int_{{\cal{K}}} u\psi \,dx \Big|_{t_1}^{t_2} 
+ \int_{t_1}^{t_2}\int_{{\cal{K}}} \Big( -u\psi_t + (\nabla u^m,\nabla\psi) - \chi u^{{\mathfrak q}-1}(\nabla v,\nabla\psi)\Big)\,dxdt = 0,
\end{equation}
\begin{equation}\label{1.4}
\int_{{\cal{K}}} v\psi \,dx \Big|_{t_1}^{t_2} 
+ \int_{t_1}^{t_2}\int_{{\cal{K}}} \Big( -v\psi_t + (\nabla v,\nabla\psi)\Big)\,dxdt
= \int_{t_1}^{t_2}\int_{{\cal{K}}} (-\alpha v + u)\psi \,dxdt.
\end{equation}

Note that the choices of $m$ and ${\mathfrak q}$ were made to ensure the integrability of the term $\chi u^{{\mathfrak q}-1}(\nabla v,\nabla u)$. In order to have sufficient integrability in the coupling term, the $u^{{\mathfrak q}-1}$ term must grow more slowly, so we need $q$ closer to 1, not larger. The threshold  $\frac {(N-2)_+}{N+2}<m<1$ also appears in other contexts. For example, the same interval appears to prove the analyticity of the solution (see \cite{DBKV}). Limiting cases can only be considered if a high integrability property is demonstrated following the approach introduced by Gianazza and Schwarzacher (\cite{GS}).

\subsection*{Context and related literature}

Chemotaxis models of porous medium type with nonlinear drift terms have attracted considerable attention in recent years (see, e.g., \cite{KL,SK} and the references therein). Such systems describe the directed movement of biological cells in response to gradients of chemical substances, known as chemoattractants. When $m=1$ and ${\mathfrak q}=2$, the system \eqref{1.1} reduces to the classical \emph{Keller--Segel model}, introduced in the seminal work \cite{KS} to describe aggregation of \emph{Dictyostelium discoideum} by chemotaxis. 
For a comprehensive survey of Keller--Segel type systems see \cite{BBTW}-
In the present case, we focus on the so-called fast diffusion regime, corresponding to diffusion of the form $D(u) = u^{m-1}$ with $0 < m < 1$.  In this regime, diffusion becomes singular near vacuum ($u = 0$), and thus enhances the spreading of cells when the local population density is low. Biologically, this can be interpreted as an adaptive mechanism by which cells avoid overcrowding and disperse more efficiently in regions of scarcity. The parameter $m$ regulates the strength of this nonlinear diffusion: smaller values of $m$ correspond to faster diffusion and a stronger regularizing effect near vanishing density. 
The drift term, on the other hand, accounts for chemotactic aggregation driven by the chemical gradient. In models of the form
\[
 u_t = \nabla \cdot (u^{m-1} \nabla u - \chi u^{q-1} \nabla v),
\]
the exponent $q$ captures nonlinear features of the chemotactic sensitivity, reflecting the fact that the cell response to chemical stimuli may saturate or amplify depending on local concentrations. The parameter $\chi>0$ measures the overall strength of chemotactic attraction. In the particular case $\chi=1$, that is our case, cells move up the gradient of the chemoattractant, corresponding to purely attractive chemotaxis with normalized sensitivity.
From a mathematical perspective, the coexistence of singular diffusion and nonlinear drift presents significant analytical challenges. The degeneracy at $u=0$ and the possible blow-up mechanisms due to aggregation compete, leading to a delicate balance between diffusion-dominated and drift-dominated regimes.  This balance strongly influences the qualitative behavior of solutions, including the existence of global in time solutions, formation of singularities, and possible pattern formation phenomena.
 Actually, if ${\mathfrak q} < 1$, the response becomes saturated (inhibition or desensitization effect). Now, for systems with singular diffusion ($m<1$), the natural tendency of the system is to strongly disperse cells in low-density regions. If very strong chemotaxis (${\mathfrak q}\gg 1$) is also introduced, competition between the two mechanisms easily leads to blow-up or non-integrability of the terms.
Consequently, the fast-diffusion chemotaxis literature generally restricts 
${\mathfrak q}$  to the range $1<{\mathfrak q}\leq 2$
 often taking ${\mathfrak q}$ close to 1 in biologically realistic settings. The hypothesis of ${\mathfrak q}>2$ is necessary in degenerate cases ($ m >1$) to balance the loss of diffusion, but in the singular case ($m < 1$) the diffusion is already strong and there is no need to amplify the aggregation.
The analysis of such systems often relies on entropy methods, energy estimates, and variational techniques adapted to the nonlinearity of the problem, extending the framework of classical Keller-Segel models to a broader class of nonlinear diffusion operators.\\
The qualitative behavior of solutions---global existence, blow-up, extinction---remains only partially understood. 
To avoid complications arising from blow-up, we restrict our analysis to a bounded time interval $(0,T]$, thereby assuming $u$ remains bounded.\\
For degenerate diffusion, there exists an extensive literature on global existence, finite-time blow-up, and H\"older regularity of weak solutions, in both parabolic--parabolic and parabolic--elliptic settings (see \cite{IY, IY3, LV, SK2, SK,W2, X}) and for the H\"older continuity to their  weak solutions, see  \cite{B, DM, 
KL, MRVVp, WMM} (see also  the references therein).\\
The singular regime has been less studied, due to the analytical difficulties created by the singularity at $u=0$. 
In the parabolic-elliptic case, Sugiyama and Yahagi \cite{SY}  assume ${\mathfrak q}$ close to 1, discuss the equilibrium between $m$ and ${\mathfrak q}$ and obtain global weak solutions for small initial data, blow-up in dimension $N=2$ with critical data, decay estimates for $1-\tfrac{2}{N} \leq m < 1$, and finite-time extinction for $0<m<1-\tfrac{2}{N}$. 
They identified the borderline case $m=1-\tfrac{2}{N}$ separating decay from extinction.  
Moreover,  Miura and Sugiyama \cite{MS}  established uniqueness of weak solutions in the H\"older class under the condition $m>\max\{ \tfrac{1}{2}-\tfrac{1}{N}, 0\}$, using a duality method combined with vanishing viscosity arguments. 
Later, Kawakami and Sugiyama \cite{KSu} extended uniqueness results for both degenerate and singular cases assuming the chemotactic term  to be quadratic or sub-quadratic.

\subsection*{Aim of the present work}
 Our results contribute to the ongoing effort to understand the fine properties of chemotaxis systems with nonlinear diffusion, and confirm that the delicate interplay between singular diffusion and aggregation admits a regularizing mechanism consistent with the porous medium paradigm.
 
In analogy with the porous medium and fast diffusion equations, where H\"older continuity is the optimal regularity for weak solutions, our goal is to establish \emph{local H\"older regularity} for weak solutions of the singular chemotaxis system \eqref{1.1} in the parabolic--parabolic case. 
Our proof combines the \emph{De Giorgi--Di Benedetto iteration method} (\cite{CDIB, degiorgi, DIB, DBGVb, DUV}) applied to the first equation in \eqref{1.1} with suitable apriori estimates for the chemoattractant $v$.\\

More generally, the first equation in \eqref{1.1} can be regarded as a prototype of
\[
u_t = \operatorname{div}(A(x,t,u,\nabla u)) - \chi \operatorname{div}(u^{{\mathfrak q}-1}\nabla v),
\]
with measurable coefficients $A(x,t,u,\nabla u)$ satisfying the standard structure conditions
\[
\begin{cases}
A(x,t,u,\nabla u)\cdot \nabla u \geq C_0 m |u|^{m-1}|\nabla u|^2, \\[6pt]
|A(x,t,u,\nabla u)| \leq C_1 m |u|^{m-1}|\nabla u|,
\end{cases}
\quad \text{a.e. in } \mathbb{R}^N \times (0,T),
\]
for some positive constants $C_0,C_1$ (see \cite{CDIB,DBGVb, DUV}). \\ The results being local, assume without loss of generality  $u\in L^{\infty}(  {\mathcal K} \times (0,T])$.
\begin{theorem}[H\"older regularity]\label{thm:holder}
Let $u$ be a locally bounded weak solution of the first equation in \eqref{1.1} under the condition \eqref{mathfrak q}. 
Then $u$ is locally H\"older continuous in $\mathbb{R}^N \times (0,T]$. 
More precisely, there exist constants $\gamma>0$ and $\alpha \in (0,1)$, depending only on the data $\{m, \; N\}$, such that for every compact set ${\mathcal{K}} \subset \mathbb{R}^N$,
\[
|u(x_1,t_1) - u(x_2,t_2)| \leq 
\gamma  \|u\|_\infty \Big( |x_1-x_2| + \|u\|_\infty^{\tfrac{m-1}{2}} |t_1-t_2|^{1/2} \Big)^{\alpha},
\]
for all $(x_1,t_1), (x_2,t_2) \in {\mathcal K} \times (0,T]$ and where $\| u \|_{\infty}$ is the $L^{\infty}$ norm on $ {\mathcal K} \times (0,T]$.
\end{theorem}

\medskip
\noindent
\textbf{Structure of the paper.}
In Section~2, we introduce the precise functional framework and recall auxiliary results from the regularity theory of degenerate and singular parabolic equations. Sections~3-6 are devoted to prove suitable energy estimates in order to adaptat  De Giorgi--Di Benedetto method (\cite{degiorgi, DIB})   to the chemotaxis system.  In Section~7 we prove logarithmic estimates. In Section~8, we prove our main H\"older continuity result.

\section{Preliminary Results}
\label{Pre}
In this section we present some auxiliary results that will be used throughout the paper. Let $K_{R}(x_{0}) \subset \mathbb{R}^{N}$ denote the $N$-dimensional cube centered at $x_{0}$ with edge length $R<1$. We write $K_{R} := K_{R}(0)$ when the cube is centered at the origin, and denote its Lebesgue measure by $|K_{R}|$. Moreover, let 
$
Q_{R}(T) := K_{R}(x_{0}) \times [0,T],
$
and, without loss of generality, assume $x_{0} = 0$.  

Consider the parabolic space $L^{q,r}(Q_{R}(T))$ endowed with the norm
\[
\|w\|_{L^{q,r}(Q_{R}(T))} := 
\left( \int_{0}^{T} \Big( \int_{K_{R}} |w|^{q} \, dx \Big)^{\frac{r}{q}} dt \right)^{\! \frac 1r} < \infty.
\]

For $p > 1$, we define
\[
V^{p}(Q_{R}(T)) := L^{\infty}(0,T; L^{p}(K_{R})) \cap L^{p}(0,T; W^{1,p}(K_{R})),
\]
with norm
\[
\|w\|_{V^{p}(Q_{R}(T))} := 
\operatorname*{ess\,sup}_{t \in (0,T)} \|w(\cdot,t)\|_{L^{p}(K_{R})} + \|\nabla w\|_{L^{p}(Q_{R}(T))} < \infty.
\]
We also define
\[
V^{p}_{0}(Q_{R}(T)) := L^{\infty}(0,T; L^{p}(K_{R})) \cap L^{p}(0,T; W^{1,p}_{0}(K_{R})).
\]

\begin{lemma}[Embedding Lemma]\label{lemma:embedding}
There exists a constant $\gamma = \gamma(N,p,s)$ such that for every function $w \in V^{p}_{0}(Q_{R}(T))$,
\begin{equation*}\label{ineq:embedding}
\iint_{Q_{R}(T)} |w(x,t)|^{q} \, dx\,dt 
\leq \gamma^{q} 
\left( \iint_{Q_{R}(T)} |\nabla w|^{p} \, dx\,dt \right)
\left( \operatorname*{ess\,sup}_{0<t<T} \int_{K_{R}} |w(x,t)|^{s} \, dx \right)^{\!\frac{p}{N}},
\end{equation*}
where $q = \frac{p(N+s)}{N}$.\\
\end{lemma}
 As a consequence for  $s=p$ Lemma \ref{lemma:embedding} takes the form
\begin{equation*}\label{ineq:Lqr}
\|w\|_{L^{q,d}(Q_{R}(T))} \leq \gamma \|w\|_{V^{p}_{0}(Q_{R}(T))},
\end{equation*}
where $q,d$ satisfy the relation
$\frac{1}{d} + \frac{N}{pq} = \frac{N}{p^{2}}.
$
In the case $1 \leq p < N$, the admissible range is
$
q \in \left[ p, \frac{Np}{N-p} \right], \ \  d \in [p, \infty].
$

\subsection*{Steklov Average}
Since the solution of the first equation in (1.1) has limited regularity with respect to the time variable, we shall make use of the \emph{Steklov average} $u^{h}$ of the weak solution $u$, defined for $h > 0$ by
\begin{equation*}\label{def:Steklov}
u^{h}(\cdot,t) := \frac{1}{h} \int_{t}^{t+h} u(\cdot,\tau) \, d\tau.
\end{equation*}
For further details we refer the reader to \cite{DIB}.  

We now recall an alternative formulation of weak solutions to (1.1): fix $t \in (0,T)$, let $h > 0$ with $0 < t < t+h < T$, and replace in \eqref{1.3} the interval $(t_{1},t_{2})$ with $(t,t+h)$. Choosing a test function $\psi$ independent of $\tau \in (t,t+h)$, dividing by $h$, and using the Steklov averages, we obtain
\begin{equation}\label{weak:Steklov}
 \int_{{\cal{K}}\times \{ t \}}  \Big[ (u^{h})_{\tau}\psi + \big( (\nabla u^{m})^{h}, \nabla \psi \big) 
- \chi \big( (u^{{\mathfrak q}-1}\nabla v)^{h}, \nabla \psi \big) \Big] \, dx = 0,
\end{equation}
for all $0<t<T-h$ and 
for all locally bounded, nonnegative test functions
\[
\psi \in W^{1,2}_{\mathrm{loc}}(0,T;L^{2}({\cal{K}})) \cap L^{2}_{\mathrm{loc}}(0,T; W^{1,2}_{0}({\cal{K}})).
\]
Fix  a subinterval $0<t_1<t_2<T$, integrating over $[t_{1},t_{2}]$ and letting $h \to 0$  \eqref{weak:Steklov} recovers \eqref{1.3}.

\begin{lemma}[De Giorgi's theoretical Lemma]\label{lemma:DeGiorgi}
Let $w \in W^{1,1}(K_{R}(x_{0}))$, and let $k < \ell$ be real numbers. There exists a constant $\gamma_{D}$ depending only on $N$, and independent of $k,\ell,w,x_{0},R$, such that
\begin{equation*}\label{ineq:DeGiorgi}
(\ell - k)\, |\{ w < k \}| \leq \gamma_{D} \frac{ R^{N+1} }{|\{ w > \ell \}| }
\int_{\{ k < w < \ell \}} |Dw| \, dx.
\end{equation*}
\end{lemma}

\begin{lemma}[Fast Geometric Convergence]\label{lemma:fastgeo}
Let $(X_{n})$ and $(Y_{n})$, $n=0,1,\dots$, be two sequences of positive numbers satisfying
\begin{equation*}\label{ineq:recurrence}
X_{n+1} \leq c \, b^{n}\big( X_{n}^{1+\alpha} + X_{n}^{\alpha} Y_{n}^{1+\kappa} \big), 
\qquad 
Y_{n+1} \leq c \, b^{n}\big( X_{n} + Y_{n}^{1+\kappa} \big),
\end{equation*}
with given constants $c,b > 1$ and $\alpha,\kappa > 0$.  
If
\[
X_{0} + Y_{0}^{1+\kappa} \leq (2c)^{-\frac{1+\kappa}{\sigma}} b^{-\frac{1+\kappa}{\sigma^{2}}},
\qquad \sigma = \min\{\kappa,\alpha\},
\]
then $(X_{n}) \to 0$ and $(Y_{n}) \to 0$ as $n \to \infty$.
\end{lemma}

\subsection*{Heat Equation with Source Term}
Consider the Cauchy problem
\begin{equation}\label{eq:heat}
\begin{cases}
v_{t} = \Delta v - av + w, & (x,t) \in \mathbb{R}^{N} \times (0,\infty), \\
v(x,0) = v_{0}(x), & x \in \mathbb{R}^{N}.
\end{cases}
\end{equation}
By classical $L^{p}$ maximal regularity results (see \cite{GT,HP,L}), we obtain:

\begin{lemma}[Heat Estimate]\label{lemma:heat}
Let $v$ be the solution of \eqref{eq:heat}.  
If $v_{0} \in W^{1,p}(\mathbb{R}^{N})$ and $w \in L^{\infty}(0,\infty; L^{p_{0}}(\mathbb{R}^{N}))$, with $1 \leq p_{0} \leq p \leq \infty$ and
\begin{equation}\label{m+1}
\frac{1}{p_{0}} - \frac{1}{p} < \frac{1}{N},
\end{equation}
then for $t \in [0,\infty)$ there exist constants $C_{0}, \overline{C}_{0} > 0$, depending only on $p,p_{0},N$, such that
\begin{equation}\label{ineq:heat}
\begin{aligned}
\|v(t)\|_{L^{p}(\mathbb{R}^{N})} &\leq \|v_{0}\|_{L^{p}(\mathbb{R}^{N})} 
+ C_{0} \|w\|_{L^{\infty}(0,T;L^{p_{0}}(\mathbb{R}^{N}))}, \\
\|\nabla v(t)\|_{L^{p}(\mathbb{R}^{N})} &\leq \|\nabla v_{0}\|_{L^{p}(\mathbb{R}^{N})} 
+ \overline{C}_{0} \|w\|_{L^{\infty}(0,T;L^{p_{0}}(\mathbb{R}^{N}))}.
\end{aligned}
\end{equation}
\end{lemma}

\begin{remark}\label{remark (k-u)_{+}}
We shall adopt the shorthand notation $(k-u) := (k-u)_{+}$ whenever $k > u$ (and analogously $u-k$ in the opposite case).
\end{remark}

\section{Energy Estimates on $\{k>u\}$} \label{Energy Estimates for $k>u$}

A fundamental tool for establishing regularity is provided by energy estimates.
We introduce the ``forward'' and ``backward'' parabolic cylinders, centered at the origin:
\[
Q^{+} := K_{R} \times [0,\theta R^{2}), 
\qquad 
Q^{-} := K_{R} \times (-\theta R^{2},0],
\]
with $0<\theta <1$. We derive the energy estimates in the backward cylinder $Q:=Q^{-}$;
the corresponding forward estimates in $Q^{+}$ follow by a standard time reversal. Although we work
with cylinders centered at $(0,0)$, the statements are translation invariant and hold for cylinders
centered at any $(x_{0},t_{0})$.\\
 We introduce the real number  \,$l,\, r > 1,\,\,0<\kappa < 1$\,\,and   $\tilde l, \tilde r$\,\, related to \,\,$l, r, \,\,\kappa$\,\, by the formulas: 
\begin{equation}\label{eq.rtil} 
 1-\frac 1{l}   = \frac{2(1+\kappa)} { \tilde l} , \ \ 1-\frac 1{r}   = \frac{2(1+\kappa) } { \tilde r}.
 \end{equation}
 Moreover for the index ${\mathfrak q}$ in the term $u^{{\mathfrak q}-1 }$ in (\ref{1.1}) we assume \eqref{mathfrak q} i.e.
 \begin{equation*}
1<{\mathfrak q}<
 \frac  {(m+1)(N+2) }{2N}. 
  \end{equation*}

\begin{lemma}[Local Energy Estimates on $\{ k>u \}$]\label{lem:energy}
Let $(u,v)$ be a locally bounded weak solution of~\emph{(1.1)}.There exist  a constant $\widehat{C}>0$ and  a constant $  I_{{\mathfrak d}}>0$  depending only on the data, such that for every cylinder $Q$,
for every piecewise smooth cut-off $\eta(x,t)$ with $0\le \eta\le 1$, $\eta|_{\partial K_{R}}=0$,
and for every level $k>0$, the following estimate holds:
\begin{align}\label{eq:energy}
& \operatorname*{ess\,sup}_{t\in(-\theta R^{2},0]}
\int_{K_{R}} \eta^{2}(k-u)^{2}\,dx
\;+\; \frac12\, m\,k^{m-1}\iint_{Q}\! \big|\nabla\!\big((k-u)\eta\big)\big|^{2}\,dx\,dt \\
& \le
\widehat{C}\, \Big\{ m\,k^{m+1}\iint_{Q}\!|\nabla\eta|^{2}\,dx\,dt
\;+\; I_{{\mathfrak d}}\,k^{2{\mathfrak q}-m-1}\!\left(
\int_{-\theta R^{2}}^{0}\! |A_{k,R}(t)|^
{   \frac{\widetilde{r}}{\widetilde{l} }  }\,dt
\right)^{ \! \frac 2{\widetilde{r}} (1+\kappa)  }\notag\\
&\qquad\quad
+\, k\!\int_{K_{R}}\!\eta^{2}(x,-\theta R^{2})(k-u)\,dx
\;+\; k^{2}\!\iint_{Q}\!|\eta\,\eta_{t}|\,dx\,dt \Big\},\notag
\end{align}
with  
\begin{align}\label{AkR}
A_{k,R}(t):=\{x\in K_{R}:\,( k^{m}-u^{m}(x,t))_{+}>0\}
\end{align}.
\end{lemma}

\begin{remark}[On the proof strategy]
The structure of the argument follows the classical pattern in the monograph
\emph{DiBenedetto--Gianazza--Vespri \cite{DBGVb}}. However, due to the significant differences introduced
by the nonlinear chemotactic drift, we present a complete proof instead of merely highlighting
the novelties. This choice preserves self-containment and readability.
\end{remark}

\begin{proof}[Proof of Lemma~\ref{lem:energy}]

We argue on the backward cylinder $Q=K_{R}\times(-\theta R^{2},0]$.\\
Let us introduce as a test function

\[
\psi = - (k^{m}-u^{m})_{+}\,\eta^{2}, \ \ for \ k>0
\]

and integrate over $\widetilde{Q}:=K_{R}\times(-\theta R^{2},t],\ -\theta R^{2}<t\le 0$.\\
By using the alternate weak formulation, modulus a standard Steklov averaging process, we obtain
\begin{align}
\label{eq:basic-identity}
&-\iint_{\widetilde{Q}} u_{\tau}(k^{m}-u^{m})\eta^{2}\,dx\,d\tau
\;+\;\iint_{\widetilde{Q}} \nabla u^{m}\cdot \nabla\big(-(k^{m}-u^{m})\eta^{2}\big)\,dx\,d\tau\\
&\;=\;\iint_{\widetilde{Q}} u^{{\mathfrak q}-1}\nabla v\cdot \nabla\big(-(k^{m}-u^{m})\eta^{2}\big)\,dx\,d\tau \notag,
\end{align}
where we used the symbol $(k^{m}-u^{m}) :=(k^{m}-u^{m})_{+}$.\\
Set $M_{1}+M_{2}=M_{3}$ for the three terms in \eqref{eq:basic-identity}.

\medskip\noindent\textbf{Step 1: Time term $M_{1}$.}
In the set $\{k>u\}$ we use the identity
\begin{equation*}\label{eq:time-primitive}
-(k^{m}-u^{m})\,u_{t} \;=\; \frac{d}{dt}\!\left(\int_{u}^{k} (k^{m}-s^{m})\,ds\right).
\end{equation*}
Hence,  integrating by parts from $-\theta R^{2}$ to $t$,
\begin{align}
\label{M_{1}}
\notag &M_{1} = \int_{K_{R}} \left[\int_{u(\cdot,t)}^{k} (k^{m}-s^{m})\,ds\right]\eta^{2}\,dx - \int_{K_{R}} \left[\int_{u(\cdot,-\theta R^{2})}^{k} (k^{m}-s^{m})\,ds\right]\eta^{2}\,dx\\
\notag&- \int_{K_{R}}\!\left(\int_{-\theta R^{2}}^{t} (\int_{u}^{k}(k^{m}-s^{m})ds ) 2\eta\eta_{\tau}d\tau\right)\!dx.
\end{align}
The standard inequalities hold (see \cite{DBGVb}, Chap.3, Prop.~9.1) :
\begin{equation}
\label{eq:key-primitive}
\frac12 m k^{m-1}(k-u)^{2}
\leq \int_{u}^{k}(k^{m}-s^{m}) ds
\leq  k^{m}(k-u)
\leq  k^{m+1}.
\end{equation}

By using \eqref{eq:key-primitive}
 we obtain
\begin{align*}
\notag &M_{1}\geq  \int_{K_{R}}\tfrac12 m k^{m-1}(k-u)^{2}\,\eta^{2}\,dx-  k^{m}\int_{K_{R}} (k-u)\,\eta^{2}(x,-\theta R^{2})\,dx\\
&-  2 k^{m}\! \iint_{\tilde Q} (k-u)\,|\eta\eta_{\tau}|\,dx\,d\tau .
\end{align*}
\medskip
\noindent\textbf{Step 2: Diffusion term $M_{2}$.}
We expand
\begin{align*}
M_{2}
= \iint_{\widetilde{Q}} \nabla u^{m}\cdot \nabla\big(u^{m}-k^{m}\big)\eta^{2} dxd\tau +2\iint_{\widetilde{Q}} \nabla u^{m}\cdot \big(u^{m}-k^{m}\big) \eta\nabla\eta dx d\tau.
\end{align*}

Using Young inequality in the second term with $\epsilon>0$, it follows
\begin{align}
\label{M2}
 M_{2}&\geq \iint_{\widetilde{Q}} |\nabla (k^{m} -u^{m})|^{2} \eta^{2} dxd\tau- \epsilon \iint_{\widetilde{Q}} |\nabla\!\big(k^m-u^m\big)|^{ 2}\eta^{2}\,dx\,d\tau \\
&-  \frac{1}{\varepsilon}\,  \iint_{\widetilde{Q}}(k^{m}-u^{m})^{2}|\nabla\eta|^{2} \notag
\end{align}
\medskip
\noindent \textbf{Step 3: Chemotaxis term $M_{3}$.}
Expand
\begin{align*}
M_{3}  &=- \iint_{\widetilde{Q}} u^{q-1}\nabla v\cdot \nabla(k^{m}-u^{m})\,\eta^{2}\,dx\,d\tau\\
&- 2\iint_{\widetilde{Q}} u^{q-1}\nabla v\cdot (k^{m}-u^{m})\,\eta\,\nabla\eta\,dx\,d\tau.
  \end{align*}

Applying the Young's inequality to the absolute values of both terms of $M_3$, for any $\delta, \hat \delta \in(0,1)$,
we get
\begin{align}
\label{M33}
M_{3}
& \leq  \frac{\delta}{2}\iint_{\widetilde{Q}}|\nabla(k^{m}-u^{m})|^{2} \eta^{2}dxd\tau
+ \frac{1}{2\delta}\iint_{\widetilde{Q}} u^{2(q-1)}\eta^{2}|\nabla v|^{2} dx d\tau  \\
&+\frac{\hat \delta}{2} \iint_{\widetilde{Q}} (k^{m}-u^{m})^{2}|\nabla\eta|^{2} dx d\tau
\;+\;  \frac{1}{2 \hat \delta}\iint_{\widetilde{Q}} u^{2({\mathfrak q}-1)}\eta^{2}|\nabla v|^{2}\,dx\,d\tau \notag 
\end{align}

\medskip\noindent\textbf{Step 4: Gathering terms and absorbing gradients.}

Now we estimate the term $\iint_{\widetilde{Q}} |\nabla (k^m-u^m)|^2 \eta^{2}\, dx\,d\tau$ in \eqref{M2} and \eqref{M33}. By using the following Caccioppoli-type bound
\begin{align*}
 |\nabla (k^m-u^m)|^2 \eta^{2}
\ge \tfrac12\,m^{2}k^{2(m-1)}|\nabla\!\big((k-u)\eta\big)|^{2} - m^{2}k^{2(m-1)}(k-u)^{2}|\nabla\eta|^{2},
\end{align*}
we can write 
\begin{align}\label{k^m-u^m eta^2}
 &\iint_{\widetilde{Q}} |\nabla (k^m-u^m)|^2 \eta^{2}\, dx\,d\tau \geq  \frac 1 2 m^{2} k^{2(m-1)} \iint_{\widetilde{Q}} |\nabla ((k-u)\eta) |^2 dx\,d\tau \\
&- m^{2}k^{2(m-1)} \iint_{\widetilde{Q}} (k-u)^{2}|\nabla \eta|^{2} dx d\tau \notag
\end{align}

Moreover,  to estimate the term $ \iint_{\widetilde{Q}} (k^{m}-u^{m})^{2}|\nabla\eta|^{2}\,dx\,d\tau$ in  \eqref{M2} and \eqref{M33}, we use \eqref{eq:key-primitive}
and take  into account of \eqref{k^m-u^m eta^2}. Then we obtain
\begin{align}\label{eq:diffusion-final}
&\int_{K_{R}} \tfrac12 m k^{m-1}(k-u)^{2}\eta^{2}(x,t)dx
+ \tfrac12 m^{2}k^{2(m-1)} \iint_{\widetilde{Q}} |\nabla\!\big((k-u)\eta\big)|^{2}dxd\tau\\
&\quad \leq C k^{2(m-1)}\iint_{\widetilde{Q}} (k-u)^{2}|\nabla\eta|^{2}dxd\tau  + C\iint_{\widetilde{Q}} u^{2({\mathfrak q}-1)}\eta^{2}|\nabla v|^{2}dxd\tau \notag\\
& +k^{m} \int_{K_{R}} (k-u)\,\eta^{2}(x,-\theta R^{2})dx + 2k^{m+1}\iint_{\widetilde{Q}} |\eta\eta_{\tau}|\,dx\,d\tau \notag.
\end{align}


\medskip\noindent\textbf{Step 5: Controlling the drift contribution.}
\\
We are left to bound
\(
\iint_{\widetilde{Q}} u^{2({\mathfrak q}-1)}\eta^{2}|\nabla v|^{2}\,dx\,d\tau.
\)
 Applying H\"older in $x$ with exponents $l$ and $l'=\frac{l}{l-1}$, 
and then H\"older in $\tau$
with exponents $r$ and $r'=\frac{r}{r-1}$
 yields
\begin{align}\label{eq:drift-raw}
&\iint_{\widetilde{Q}} \eta^{2}u^{2({\mathfrak q}-1)}|\nabla v|^{2}\,dx\,d\tau \notag\\
&\le k^{2({\mathfrak q}-1)}
\left(\int_{-\theta R^{2}}^{t}
\Big(\int_{K_{R}} |\nabla v|^{2l}\,dx\Big)^{\! \frac rl} d\tau \right)^{\! \frac 1r}
\left(\int_{-\theta R^{2}}^{t} |A_{k,R}(\tau)|^{\frac {(l-1)r}{l(r-1)}}\,d\tau\right)^{\!\frac {r-1}r}\notag\\
&\le k^{2({\mathfrak q}-1)}\, I_{{\mathfrak d}}\,(\theta R^{2})^{\frac 1r}
\left(\int_{-\theta R^{2}}^{0} |A_{k,R}(t)|^{\frac {(l-1)r}{l(r-1)}}\,dt\right)^{\! \frac{r-1}r},
\end{align}
where, by using \eqref{m+1}, \eqref{ineq:heat}  in  Lemma \ref{lemma:heat}
with $p_0=m+1, \ p=2l$ and
 $ \frac{1}{m+1}- \frac{1}{2l}< \frac{1}{N},\,$
 the data-dependent constant $I_{{\mathfrak d}}$ is
\[
I_{{\mathfrak d}}:=
\left[
C_{0}\,\sup_{-\theta R^{2}\le t\le 0}\!\left(\int_{K_{R}} |u(x,t)|^{m+1}\,dx\right)^{\frac 1{m+1}}
\;+\;\left(\int_{K_{R}} |\nabla v_{0}(x)|^{2l}\,dx\right)^{\!\frac 1{2l}}
\right]^{\!2}.
\]
We recall that $A_{k,R}(t)$ present in \eqref{eq:drift-raw} is the active set  defined in \eqref{AkR}.\\
Rewriting \eqref{eq:drift-raw} with the exponents $\widetilde{r},\widetilde{l}$ (see \eqref{eq.rtil}),  observing  that
  $(\theta R^{2})^{\frac 1r}=(\theta R^2)^ {1-\frac{2}{\tilde{r}}(1+\kappa)} <1,$ we obtain \begin{equation}\label{eq:drift-final}
\iint_{\widetilde{Q}} \eta^{2}u^{2({\mathfrak q}-1)}|\nabla v|^{2}\,dx\,d\tau
\;\le\; I_{\mathfrak d} \,k^{2({\mathfrak q}-1)}\,\!
\left(\int_{-\theta R^{2}}^{0}\! |A_{k,R}(t)|^{\frac{\widetilde{r}} {\widetilde{l}}}\,dt\right)^{\! \frac 2{\tilde r}(1+k)},
\end{equation}

(which is the scaling used in the
iteration), this yields exactly the second right-hand side term in \eqref{eq:energy} after
dividing the whole inequality by $m k^{m-1}$.

\medskip\noindent\textbf{Step 6: Conclusion.}
Take the essential supremum in time of the leftmost term in \eqref{eq:diffusion-final},
replace $\widetilde{Q}$ with $Q$, divide the resulting inequality by $m k^{m-1}$,
and use \eqref{eq:diffusion-final} together with \eqref{eq:drift-final}.
This yields \eqref{eq:energy} with a possibly different constant $\widehat{C}>0$
depending only on the structural data (and on the choice of the universal parameters
in the Young inequalities). The proof is complete.
\end{proof}

\section{A De~Giorgi Type Lemma on $\{k > u\}$}

For a cylinder $Q_{R}(\theta)$, we denote
\[
\mu^{+} := \operatorname*{ess\,sup}_{Q_{2R}(\theta)} u, 
\qquad 
\mu^{-} := \operatorname*{ess\,inf}_{Q_{2R}(\theta)} u,
\qquad 
\omega := \mu^{+} - \mu^{-}.
\]
Let $\xi, a \in (0,1)$ be fixed numbers and assume $R \leq \xi$.

\begin{lemma}[De~Giorgi type lemma on $\{k > u\}$]\label{lemma:degiorgi}
There exists a positive number $\nu \in (0,1)$, depending only on $\theta, \omega, \xi, a$ and the data, such that if
\begin{equation}
\label{Firstip} 
\bigl|\{(x,t)\in Q_{2R}(\theta) : u < \mu^{-} + \xi\omega\}\bigr| \leq \nu\, |Q_{2R}(\theta)|,
\end{equation}
then
\begin{equation*}
\label{First} 
u > \mu^{-} + a\xi \omega \quad \text{a.e.\ in } Q_{R}(\theta).
\end{equation*}
\end{lemma}

For the reader's convenience we provide a detailed proof.

\begin{proof}
For $n=0,1,2,\dots$, set
\[
R_{n} := R + \frac{R}{2^{n}}, 
\qquad K_{n} := K_{R_{n}}, 
\qquad Q_{n} := K_{n} \times (-\theta R_{n}^{2},0],
\]
all centered at the origin $(0,0)$.  
The levels are chosen as
\[
k_{n} := \mu^{-} + \xi_{n}\omega, 
\qquad 
\xi_{n} := a\xi + \frac{1-a}{2^{n}}\xi.
\]
To make the proof simpler, we assume $\mu^-=0$.\\
The cutoff functions $\eta_{n}(x,t) = \eta_{1,n}(x)\eta_{2,n}(t)$ are defined as follows:
\begin{itemize}
  \item $\eta_{1,n} = 1$ in $K_{n+1}$, $\eta_{1,n} = 0$ in $K_{n}\setminus K_{n+1}$, with 
  $|\nabla \eta_{1,n}| \leq \frac{2^{n+1}}{R}$;
  
  \item $\eta_{2,n} = 0$ for $t < -\theta R_{n}^{2}$, $\eta_{2,n} = 1$ for $t \geq -\theta R_{n+1}^{2}$, and $
  0 \leq (\eta_{2,n})_{t} \leq \frac{2^{2(n+1)}}{\theta R^{2}}.
  $
\end{itemize}

Applying inequality \eqref {eq:energy}     over $Q_{n}$ and defining
\begin{equation*}
\bigl| A_{n} \bigr| := \bigl|\{u < k_{n}\}\cap Q_{n}\bigr|
= \int_{-\theta R_{n}^{2}}^{0} \bigl|A_{n}(t)\bigr|\, dt, \quad 
  A_{n}(t) := \{x \in K_{R_{n}} : u(x,t) < k_{n}\},
\end{equation*}
we obtain (after suitable estimates involving the cutoff functions and the bound $\xi\omega < 1$) the inequality
\begin{align}
\label{eq:EE4}
&\operatorname*{ess\,sup}_{-\theta R_{n}^{2} < t \leq 0} 
\int_{K_{R_{n}}} \eta_{n}^{2} (k_{n}-u)^{2}\, dx
+ \frac{1}{2} m(\xi \omega)^{m-1} 
\iint_{Q_{n}} |\nabla ((k_{n}-u)\eta_{2,n})|^{2}\, dx\,dt \\
&  \leq 
c_{1} \frac{36m^{2} 2^{2(n+1)}}{R_{n}^{2}} (\xi\omega)^{m+1}
\Biggl(1 + \frac{1}{\theta(\xi\omega)^{m-1}} \! \Biggr) \! \Biggr(   |A_{n}|\!
+ \! \Bigl( \int_{-\theta R_{n}^{2}}^{0} |A_{k_{n},R_{n}}|(t)^{\frac{\tilde r} {\tilde l}}\, dt \Bigr)^{\frac 2{\tilde r} (1+\kappa)} \! \Biggr). \nonumber
\end{align}

Moreover, in $Q_{n+1}$ we have
\[
|(k_{n}-u)|^{2} \geq |k_{n}-k_{n+1}|^{2} = \bigl(\xi\omega (1-a)2^{-(n+1)}\bigr)^{2}.
\label{eq:EE5}
\]
Therefore,
\[
\iint_{Q_{n}} |(k_{n}-u)\eta_{n}|^{2}\, dx\,dt 
\geq \Bigl( (1-a)\xi\omega 2^{-(n+1)} \Bigr)^{2} |A_{n+1}|.
\label{eq:EE6}
\]

Using H\"older's inequality together with Lemma \ref{lemma:embedding}(with $p=2$, $q=\tfrac{2N+2}{N}$) and recalling that $\eta_{n}=1$ in $Q_{n+1}$, we arrive at
\begin{align}
\label{eq:EE7}
&\Bigl((1-a)\xi\omega 2^{-(n+1)}\Bigr)^{2}|A_{n+1}|
\leq \iint_{Q_{n+1}\cap \{u<k_{n+1}\}} (k_{n}-u)^{2}\, dx\,dt  \\
&\leq \gamma \!
\Biggl( \iint_{Q_{n}}\! \! |\nabla((k_{n}-u)\eta_{n})|^{2}\, dx\,dt \Biggr)^{\frac N{N+2}} \!\!
\Biggl( \operatorname*{ess\,sup}_{-\theta R_{n}^{2} \leq t \leq 0} 
\int_{K_{R_{n}}} \! \! \eta_{n}^{2}(k_{n}-u)^{2}\, dx \Biggr)^{\frac 2{N+2}}
|A_{n}|^{\frac 2{N+2}}.\nonumber
\end{align}
Let us introduce the sequences
 \begin{equation*}
 X_n= \frac {| A_n|} {|Q_n| } \quad  \text{and} 
 \quad Y_n=
   \frac{ \Bigg( \int^0_{- \theta R_n^2}  | A_n(t)|^{\frac {\tilde r} {\tilde l}}dt \Bigg)^{\frac {2 } {\tilde r }  }} {|Q_n|^{\frac{N}{N+2}}} .
\end{equation*}

Combining estimates \eqref{eq:EE4}-\eqref{eq:EE7}, and iterating the resulting recursive inequalities for the normalized sequence we get
\begin{equation}
\label{EE14}      
\begin{aligned}
&X_{n+1}
\leq C_1 \ 16^{n+1}\Big(X_n^{1+\frac{2}{N+2}}+ X_n^{\frac{2}{N+2}}Y_n^{1+\frac{2}{N}} \Big)
\end{aligned}
\end{equation}
with $ C_1>0.$
Now we    derive an  estimate of the same type  for $Y_n$. Starting from the inequality
\begin{equation*}
\begin{aligned}
 & Y_{n+1}(k_{n}- k_{n+1} )^2 \leq\frac{2^N}{\theta^{\frac{N}{N+2}}|K_{R_{n}}|}\Bigg[ \int_{- \theta R_{n+1}^2}^0  \! \Bigg( \int_{K_{R_{n+1}}}(\eta_n(k_n-u))^{\tilde{l}}\,\, dx\Bigg)^{\frac{\tilde{r}}{\tilde{l}}} dt \Bigg]^{\frac{2}{\tilde{r}}}, 
 \end{aligned}
\end{equation*}
applying Lemma \ref{lemma:embedding},
 after some calculations we arrive to 
 \begin{equation}
\label{2Yn}      
\begin{aligned} Y_{n+1}\leq C_2 16^{n+1} \Big(X_n  + Y_n^{1+\frac{2}{N}}\Big).   
 \end{aligned}
\end{equation}

From  \eqref{EE14} and   \eqref{2Yn}  by   Lemma \ref{lemma:fastgeo},  
it follows that $X_{n}\to 0$ and $Y_{n}\to 0$ as $n\to\infty$, provided the initial smallness condition
\[
X_{0} + Y_{0}^{1+\frac 2 N} \leq
(2C_3)^{-\frac{1+\frac{2}{N}}{\sigma}}16^{-\frac{1+\frac{2}{N}}{\sigma^2}} := \nu, \ \ \sigma := \tfrac{2}{N+2} \ and \  \kappa= \frac 2 N
\]
is satisfied with $C_3=\max\{C_1,C_2\}.$
Hence, inequality \eqref{Firstip}  follows, and the lemma is proved.
\end{proof}

\begin{remark}
Lemma~\ref{lemma:degiorgi} represents a key step in the De~Giorgi iteration scheme.  
It shows that, under a quantitative smallness assumption on the set where the solution remains close to its essential infimum, one can improve the lower bound of the solution in a smaller cylinder.  
Iterating this procedure eventually leads to the reduction of the oscillation of $u$ inside shrinking cylinders, which is the cornerstone of proving local boundedness and, subsequently, H\"older continuity.  
This mechanism, originally devised by De~Giorgi for elliptic equations, extends here to the present parabolic setting with nonlinear diffusion and drift.
\end{remark}

\section{Energy Estimates on $\{u> k\}$}

As in Section \ref{Energy Estimates for $k>u$}, consider the cylinder 
\[
Q = Q_{R}(\theta) := K_{R} \times (-\theta R^{2}, 0].
\]

\begin{lemma}[Local Energy Estimates on $\{u> k\}$]\label{lem:energy2}
Let $(u,v)$ be a locally bounded weak solution of problem \eqref{1.1}. 
  There exist a  constant $\gamma>0$ and  a constant $I_{{\mathfrak d}}>0$ depending only on the data, such that, for any cut-off function $\eta \geq 0$ vanishing on $\partial K_{R}$, for every level $k > 0$ and for every cylinder $Q$, the following inequality holds:
\begin{align}
\label{EEu>k}
& \underset{-\theta R^2<t\leq 0}{\text{ess\,sup}} \int_{K_{R}} (u-k)^{2}\eta^{2}(x,t)\,dx 
+ \iint_{Q} u^{m-1}\,|\nabla((u-k)\eta)|^{2}\,dx\,dt  \\
&  \leq \gamma \iint_{Q} u^{m-1}(u-k)^{2}|\nabla \eta|^{2}\,dx\,dt 
+ I_{{\mathfrak d}}\,(\mu^+)^{\,2{\mathfrak q}-m-1}   \left(\int_{-\theta R^{2}}^{0} 
\big|A_{k,R}(t)\big|^{\frac{\tilde r}{ \tilde \ell}}\,dt \right)^
{\!\frac{2}{\tilde r}(1+\kappa)} \nonumber   \\
& \qquad \quad  + \int_{K_{R}} (u-k)^{2}\eta^{2}(x,-\theta R^{2})\,dx 
+ \iint_{Q} (u-k)^{2}\eta\,\eta_{t}\,dx\,dt,  \nonumber
\end{align}
\end{lemma}
with $A_{k,R}(t)=\{x\in K_{R}:\,(u^{m}-k^{m}(x,t))_+>0\} $.
\begin{proof}
For $u>k$, consider the test function 
$
\psi = (u^m-k^m)\eta^{2}
$
over the truncated cylinder $\widetilde{Q} = K_{R}\times(-\theta R^{2},t]$, with $t \leq 0$.  
From the definition of weak solution \eqref{1.3}, integrating in time from $-\theta R^{2}$ to $t$ and letting $h\to 0$ in the Steklov averages, we obtain
\begin{align}
&\iint_{\widetilde{Q}} (u^m-k^m)\eta^{2} u_{t}\,dx\,d\tau 
+ \iint_{\widetilde{Q}} \nabla u^{m}\cdot \nabla((u^m-k^m)\eta^{2})\,dx\,d\tau \nonumber \\
&= \iint_{\widetilde{Q}} u^{q-1}\nabla v \cdot \nabla((u^m-k^m)\eta^{2})\,dx\,d\tau.
\notag
\end{align}

In our proof we follow the steps in Lemma \ref{lem:energy}  with the  difference that here $u$ is estimated from above  by $\mu^{+}$, while there  (in Lemma  \ref{lem:energy}) $u$ is estimated by  $k$.
\end{proof}

\section{A De Giorgi-Type Lemma on $\{u> k\}$}

Let $\xi$ and $a$ be fixed numbers in $(0,1)$.
\begin{lemma}[De Giorgi-Type Lemma for $u > k$]
\label{DGbisLemma}
Assume that
\begin{equation}\label{eq:mu_condition}
\mu^+ \leq \frac{13}{12}\,\omega.
\end{equation}
There exists a positive constant $\nu^*$, depending only on $\theta, \omega, \xi, a$, and the data, such that if
\begin{equation*}\label{eq:measure_condition}
\left| \{ (x,t) \in Q_{2R}(\theta) : u > \mu^+ - \xi \omega \} \right| \leq \nu^* |Q_{2R}(\theta)|,
\end{equation*}
then
\begin{equation}\label{eq:degiorgi_conclusion}
u \leq \mu^+ - a \xi \omega \quad \text{a.e. in } Q_R(\theta).
\end{equation}
\end{lemma}

\begin{proof}
As in Lemma \ref{lemma:degiorgi}, let us introduce the sequences of cubes $(K_n) = (K_{R_n})$, cylinders $(Q_n)$, and cut-off functions $\eta_n = \eta_{1,n}(x)\eta_{2,n}(t)$.  
Define the sequence of truncation levels
\[
k_n = \mu^+ - \xi_n \omega, 
\qquad 
\xi_n = a \xi + \frac{1-a}{2^n}\xi.
\]
We apply the Energy Estimates from Lemma \ref{lem:energy2}
 over the cylinders $Q_n$ for the truncated function $(u-k_n)$, where
 the function $u$ is estimate from above by $\mu^+$. Moreover from the assumption 
 \eqref{eq:mu_condition}, since $\omega<1$, then $u\leq \frac {13}{12}$.  
 So we obtain
\begin{align}
 \label{eq:EE} 
 &\underset{( -\theta R^2, 0]}   {ess\,sup } \   \int_{K_{R}} (u-k)^2 \eta^2(x,t)dx +c  \Big(\frac {13}{ 12}\Big)^{m-1}  \iint_{\widetilde{Q}} |\nabla (u-k)|^2     \eta^{2}\,dx\,d\tau \\
&\leq \gamma \  \Big[ \Big(\frac {13}{ 12}\Big)^{m+1}  \iint_{\widetilde{Q}} |\nabla\eta|^{2}\,dx\,d\tau +
\Big(\frac {13}{ 12}\Big)^{2 \mathfrak{q}-m-1}
\iint_{\widetilde{Q}} \eta^{2}|\nabla v|^{2}\,dx\,d\tau \nonumber \\ 
&+   \frac {13}{ 12}  \int_{K_{R}} (u-k)   \eta^2(x, -\theta R^2)dx+ 
 \Big( \frac {13}{ 12}\Big)^{2} \iint_{\widetilde{Q}} (u-k)   \eta |\eta_{\tau}|\,dx\,d\tau\Big]  \nonumber
\end{align}
 We rewrite \eqref{eq:EE} as
\begin{align}
\label{eq:energy_Qn}
&\underset{-\theta R_n^2<t\leq 0}{\text{ess\,sup}}  \int_{K_n} (u-k_n)^2 \eta_n^2(x,t)\,dx 
+ c_1 \iint_{Q_n}  \, |\nabla((u-k_n)|^2  \eta_n^2 \,dx\,dt  \\
&\qquad \leq \gamma \left( P_1 + P_2 + P_3 + P_4 \right).\nonumber
\end{align}
with $c_1=  2 c  (\frac {13}{12})^{m-1} $. We will disgard $c_1$ .
Here $
P_3 = 0 \quad \text{(by the structure conditions)},
$
\[
P_1 \leq  \tilde c_1  \frac{2^{2(n+1)}}{ R_n^2}  
\iint_{Q_n} \chi_{[u>k_n]} \, dx dt \leq \tilde c_1  \frac{2^{2(n+1)}}{\theta R_n^2}  |A_n|
\]
with $\tilde c_1=2 \gamma  (\frac {13}{ 12} )^{m+1}$,  with   $\theta <1 $   in   the  numerator\\
\[
P_2 \leq  2 \ I_{\mathfrak d} \left(\tfrac{13}{12}\right)^{2{\mathfrak q}-m-1} 
\left( \int_{-\theta R_n^2}^{0} |A_{k_n,R_n}(t)|^{\frac {\tilde r}{\tilde l} }\,dt \right)^{\frac 2{\tilde r} (1+\kappa)},
\]
while
\[
P_4 \leq   2  \left(\tfrac{13}{12}\right)
  \frac{2^{2(n+1)}}{\theta R_n^2} (\xi \omega) 
\iint_{Q_n} \chi_{[u > k_n]} \, dx dt \leq \tilde c_2   \frac{2^{2(n+1)}}{\theta R_n^2} |A_n|,
\]
with $\xi,\omega<1.$ 
\ Also for a new constat $\gamma_3$ we get from \eqref{eq:energy_Qn}

\begin{equation}
\label{DG11}
\begin{aligned}
&\underset{-\theta R_n^2<t\leq 0}{\text{ess\,sup}}  \int_{K_n} (u-k_n)^2 \eta_n^2(x,t)\,dx 
+  \iint_{Q_n}  \, |\nabla((u-k_n)|^2  \eta_n^2 \,dx\,dt \nonumber \\
&\qquad \leq \gamma_3\Big[  \frac{2^{2(n+1)}}{\theta R_n^2}  |A_n|+ 
\left( \int_{-\theta R_n^2}^{0} |A_{k_n,R_n}(t)|^{\frac {\tilde r}{\tilde l} }\,dt \right)^{\frac 2{\tilde r} (1+\kappa)}\Big]
\end{aligned}
\end{equation}

Now, since in $Q_{n+1}$ we have
\begin{equation*}
|u-k_n|^2 \geq |k_n-k_{n+1}|^2, 
\qquad |k_n-k_{n+1}| = \xi\omega (1-a)\, 2^{-(n+1)},
\end{equation*}
it follows that
\begin{equation*}
\iint_{Q_{n+1}} (u-k_n)^2 \,dx\,dt \geq \left( (1-a)\xi\omega\,2^{-(n+1)} \right)^2 |A_{n+1}|.
\end{equation*}

Apply on the left of the last inequality first  H\"older inequality  and  then  Lemma \ref{lemma:embedding} with $ p=2,\, \ s=2, \ q = 2\frac{N+2}{N} $ to get
\begin{equation}
\label{EE8bis}
\begin{aligned}
&\Big(\frac {(1-a) \xi \omega } { 2^{n+1}}\Big)^2|A_{n+1}| 
 \leq {\gamma}^{\frac{N}{N+2}  } \Big( \iint_{Q_n}  |\nabla (u-k_n)\eta_n|^2dx dt\Big)^{\frac{N}{N+2}} \\
&\times   \underset{-\theta R_n^2 <  t \leq 0}{\text ess\, sup}\int_{K_{n}}  (u-k_n)^2 \eta_n^2 dx \Big)^{\frac{2}{N+2}} 
\   |A_n|^{ \frac {2}{N+2}}.
\end{aligned} 
\end{equation}
 
From \eqref{EE8bis} we  derive for a suitable $\gamma_5>0$, with  $ \omega,\theta,\xi \leq 1,$
\begin{equation*}
\label{DG6}
\begin{aligned}
&
|A_{n+1}|
\leq  \gamma_5 \frac {2^{2(n+1)}} {(1-a)^2  }
 \Bigg\{\frac{1}{\theta R_n^2} |A_n|^{1+ \frac {2}{N+2}}
  + \frac 1{ \theta R_n^2} \Big(\int_{-\theta R_n^2}^0 |A_{k_n,R_n}|^{\frac{\tilde{r}}{\tilde{l}}}dt|\Big)^{\frac{2}{\tilde{r}}(1+\kappa)}|A_n|^{\frac{2}{N+2}}\Bigg\}.
\end{aligned}
\end{equation*}
Now, as in Lemma \ref{lem:energy}, we introduce  the sequences
\vspace{-0.8em}
 \begin{equation*}
 \label{2Xn}
 X_n= \frac {| A_n|} {|Q_n| } \quad  and \quad
  Y_n=
   \frac{ \Bigg( \int^0_{- \theta R_n^2}  | A_n(t)|^{\frac {\tilde r} {\tilde l}}dt \Bigg)^{\frac {2 } {\tilde r }  }} {|Q_n|^{\frac{N}{N+2}}} .
\end{equation*}
which   satisfy the nonlinear iteration scheme of Lemma \ref{lemma:fastgeo} with appropriate parameters. Hence, provided the smallness condition
\vspace{-0.5em}
\begin{equation*}
\label{X0Y0}
X_0 + Y_0^{1+\tfrac{2}{N}} \leq \nu^*
\end{equation*}
holds for some $\nu^* > 0$ depending only on $\theta, \omega, \xi, a$ and the data, we conclude that
\[
X_n,\, Y_n \to 0 \quad \text{as } n\to \infty,
\]
which implies \eqref{eq:degiorgi_conclusion}.  
\end{proof}


\section{Local Logarithmic Estimates on $\{u > k\}$}

In order to derive further estimates, we introduce the logarithmic function
\begin{equation*}\label{eq:psi}
    \psi(H) = \psi(u) := \log^{+}\!\left( \frac{H}{H - (u-k) + c} \right), 
\end{equation*}
$ H = \underset{Q}{\text{ess\,sup}} (u-k),$
 $0 < c < \min\{1,H\}$ and $\log^+\! s := \max\{\log s, 0\}$ for $s > 0$.

Let $\zeta$ be a nonnegative, piecewise smooth cutoff function, independent of $t$.  
We introduce the following test function in the weak formulation of \eqref{1.1}:
\begin{equation*}\label{eq:test}
    \varphi(u) = \frac{\partial}{\partial u}\!\left( \psi^{2}(u) \right)\zeta^{2} 
    = 2\psi \psi' \zeta^{2}.
\end{equation*}
The following estimates are readily verified:
\begin{equation}\label{eq:psi-bounds}
    \psi \leq \log \frac{H}{c}, 
    \qquad \psi' \leq \frac{1}{c}.
\end{equation}

\begin{lemma}\label{lem:log-estimate}
There exists a positive constant $\gamma$, depending only on the data, such that for every cylinder $Q$ and every level $k \geq 0$, one has
\begin{align}\label{eq:lemma}
&  \sup_{-\theta R^{2} < t \leq 0}\int_{K_{R}} \psi^{2}(x,t)\zeta^{2}\, dx 
    \leq \int_{K_{R}} \psi^{2}(x,-\theta R^{2})\zeta^{2}\, dx  \\
&    \quad + \gamma \left( 
        \frac{1}{c^{2}}\!\left[ 1 + \log \frac{H}{c}\right] 
        I_{{\mathfrak d}}\, \mu^{2{\mathfrak q}-m-1} \  |K_{R}| + \log\frac{H}{c}\iint_{Q} u^{m-1}|\nabla\zeta|^{2}\, dxdt
    \right). \notag
\end{align}
\end{lemma}

\begin{proof}
Consider the cylinder $\widetilde{Q} = K_{R}\times(-\theta R^{2}, t]$, with $-\theta R^{2}< t \leq 0$.  
From the weak formulation of \eqref{1.1}, we obtain
\begin{align*}
 &\iint_{\widetilde{Q}} u_{t}(2\psi\psi')\zeta^{2}\, dxd\tau
 = - \iint_{\widetilde{Q}} \nabla u^{m}\cdot\nabla(2\psi\psi'\zeta^{2})\, dxd\tau\\
 &  + \iint_{\widetilde{Q}} u^{q-1}\nabla v\cdot\nabla(2\psi\psi'\zeta^{2})\, dxd\tau 
  =: -I_{1} + I_{2}. \notag
\end{align*}

\paragraph{Estimate of $I_{1}$.}
Expanding the gradient term, we find
\begin{align*}
-I_{1} 
&= -m \iint_{\widetilde{Q}} u^{m-1}\nabla u 
    \left( \psi'^{2}\zeta^{2}\nabla u + 2\psi\psi''\zeta^{2}\nabla u + 4\psi\psi'\zeta\nabla\zeta \right) \, dxd\tau  \notag \\
&\leq -2m\!\iint_{\widetilde{Q}} u^{m-1}\psi'^{2}(1+\psi)\zeta^{2}|\nabla u|^{2}\, dxd\tau  
  + 4m\!\iint_{\widetilde{Q}} u^{m-1}\psi\psi'\zeta|\nabla u||\nabla \zeta|\, dxd\tau \notag \\
&\leq -2m\!\iint_{\widetilde{Q}} u^{m-1}\psi'^{2}(1+\psi)\zeta^{2}|\nabla u|^{2}\, dxd\tau  \notag\\
&+\frac{m}{2}\!\iint_{\widetilde{Q}} u^{m-1}\psi\psi'^{2}\zeta^{2}|\nabla u|^{2}\, dxd\tau
    + 8m\!\iint_{\widetilde{Q}} u^{m-1}\psi |\nabla\zeta|^{2}\, dxd\tau, 
\end{align*}
where we used the Young's inequality.

\paragraph{Estimate of $I_{2}$.}
Similarly, we have
\begin{align*}
    I_{2} 
    &\leq \iint_{\widetilde{Q}} \Big(u^{{\mathfrak q}-1}\nabla v \cdot 
        \big(2\psi'^{2}(1+\psi)\zeta^{2}\nabla u\Big) + \Big(4   u^{{\mathfrak q}-1}\nabla v     \psi\psi'\zeta\nabla\zeta \Big) dxd\tau  \notag \\
  &\leq \iint_{\widetilde{Q}} u^{m-1}(1+\psi)\psi'^{2}\zeta^{2}
        \left( \tfrac{16}{m}u^{2{\mathfrak q}-2m}|\nabla v|^{2} + \tfrac{m}{16}|\nabla u|^{2}\right)\, dxd\tau  \notag\\
        & +\iint_{\widetilde{Q}}\left( 
        \tfrac{32}{m}u^{2{\mathfrak q}-m-1}\zeta^{2}(1+\psi)\psi'^{2}|\nabla v|^{2} 
        + \tfrac{m}{8}u^{m-1}\psi|\nabla\zeta|^{2}\right)\, dxd\tau.  
\end{align*}

\paragraph{Combination of terms.}
The terms involving $|\nabla u|^{2}$ combine into a nonpositive contribution and can be neglected.  

Now for the terms involving $|\nabla v|^{2}$ we use the estimates from \eqref{eq:drift-final}, with our choice of $\tilde l,\tilde r, m, \kappa,$ 
then we add the terms involving 
 $|\nabla\zeta|^{2}$.

\paragraph{Conclusion.}

The claimed inequality \eqref{eq:lemma} follows with a constant $\gamma>0$ depending only on the data.
\end{proof}


\section{H\"older Continuity}
\subsection{Preliminaries for the Proof of Theorem~\ref{thm:holder}}

The key idea in proving H\"older continuity is an \emph{induction argument}.  

Fix a point $(x_0,t_0)$ and parameters $\delta \in (0,1)$, $b>1$, $R>0$, and $\omega>0$.  
We define the sequences
\begin{equation}
\label{omegan}
R_n = \frac{R}{b^n}, \quad R_0 = R, 
\qquad \omega_{n+1} = \delta \omega_n, \quad \omega_0 = \omega,
\quad n=0,1,2,\dots
\end{equation}

and the cylinders
\[
Q_n := Q_{R_n} = K_{R_n}(x_0) \times \bigl(t_0 - \omega_n^{1-m} R_n^2, \; t_0 \bigr).
\]
\begin{proposition}
\label{prop}
If $\delta$ and $b$ can be chosen independently of $u$ and $(x_0,t_0)$, depending only on the data, such that

\begin{equation}
\label{Qn }
Q_{n+1} \subset Q_n, \qquad 
\operatorname*{ess\,osc}_{Q_n} u \leq \omega_n,
\end{equation}

then $u$ is H\"older continuous at $(x_0,t_0)$. Consequently, Theorem~\ref{thm:holder} holds (see \cite{DIB}, \cite{DBGVb}).
\end{proposition}

Assume $(x_0,t_0)=(0,0)$ and consider
$
Q_0 = K_{R_0} \times (-R_0^2,0].
$
Define
\begin{equation*}
\label{mu0+mu0- }
\mu_0^+ = \operatorname*{ess\,sup}_{Q_0} u, \qquad 
\mu_0^- = \operatorname*{ess\,inf}_{Q_0} u, \qquad 
\omega_0 = \mu_0^+ - \mu_0^- = \operatorname*{ess\,osc}_{Q_0} u. \end{equation*}

If $\omega_0 \leq 1$, then \eqref{Qn }  trivially holds for $n=0$.  
We proceed by induction. Without loss of generality, assume
\begin{equation}
\label{lambda}
R < \frac{\omega^2}{\lambda}, \quad \lambda > 1. 
\end{equation}
If \eqref{lambda} fails, then $R \geq \omega^2/\lambda$ and the oscillation is already comparable to $R$, so there is nothing to prove.
Suppose  \eqref{Qn } holds for some $n$. In the cylinder $Q_n = Q_{R_n}(\theta)$, define
\begin{equation*}
\label{mu+mu-bis}
\mu_n^+ = \operatorname*{ess\,sup}_{Q_{R_n}(\theta)} u, 
\qquad \quad 
\mu_n^- = \operatorname*{ess\,inf}_{Q_{R_n}(\theta)} u. 
\end{equation*}

We consider two alternatives: 

\paragraph{First Alternative.} By Lemma \ref{lemma:degiorgi} (with $\xi=a=\tfrac{1}{2}$), there exists $\nu>0$ such that if
\begin{equation}
\label{umagg}
\bigl|\{u \leq \mu_n^- + \tfrac{1}{2}\omega_n\} \cap Q_{\frac{R_n}2}(\theta)\bigr|
\leq \nu \ Q_{\frac{R_n}2}(\theta)|, 
\end{equation}
then
\begin{equation}
\label{umaggg}
u \geq \mu_n^- + \tfrac{1}{4}\omega_n \quad \text{a.e. in } Q_{\frac {R_n}4}(\theta). 
\end{equation}

\paragraph{Second Alternative.} If \eqref{umagg} fails, then we may assume
\begin{equation}
\frac{\omega_n}{2} < \mu_n^+ - \frac{\omega_n}{4} < \frac{5}{6}\omega_n. \label{tr34}
\end{equation}
If the left inequality fails, then  by \eqref{omegan}
\[
\operatorname*{ess\,osc}_{Q_{n+1}} u \leq \tfrac{3}{4}\omega_n = \tfrac{3}{4\delta}\omega_{n+1},
\]
and if $\tfrac{3}{4\delta}<1$, the induction closes.  

Following \cite{DIB} Ch.~III, Sec.~7 and \cite{DBGVb} Appendix B, Lemma B.9.2, there exists a time $s$ with
\[
- \theta_n \bigl(\tfrac{R_n}{2}\bigr)^2 \leq s \leq - \tfrac{1}{2}\nu \theta_n \bigl(\tfrac{R_n}{2}\bigr)^2
\]
such that
 \begin{equation}
\label{B.9.5}
|\{u(.,s)< \mu_n^-+\frac{1}{2}\omega_n\}\cap K_{\frac{R_n}2}|>\frac{1}{2}\nu |K_{\frac{R_n}2}|.
\end{equation}
This implies
\begin{equation}
\label{B.9.6}
|\{u(.,s)>\mu_n^+-\frac{1}{4}\omega_n\}\cap K_{\frac{R_n}2}|\leq(1-\frac{1}{2}\nu) |K_{\frac{R_n}2}|.
\end{equation}
From here one propagates \eqref{B.9.6} forward using logarithmic estimates (Lemma \ref{lem:log-estimate}), leading to:

\begin{lemma}
\label{n*nu}
There exists an integer $n_*$, depending only on the data and $\nu$, such that $\forall\, s<t<0$
 \begin{equation}
\label{nu*}
\bigl|\{u(\cdot,t) > \mu_n^+ - \tfrac{\omega_n}{2^{n_*}}\} \cap K_{\frac {R_n}2}\bigr|
\leq \Bigl(1 - \tfrac{\nu^2}{4}\Bigr)|K_{\frac {R_n}2}|. 
\end{equation}

\end{lemma}

\begin{proof}
  Let us consider the logarithmic estimate  \eqref{eq:lemma} for the function  $(u-k_n) $   
applied to
   the cylinder $K_{\frac{R_n}{2}}\times(s,0)$\, and the level $ k_n = \mu_n^+ -\frac{1}{4}\omega_n,$
   
 and   $c= \frac{\omega_n}{2^{i+2}},$ where $i$ will  be chosen later on.
 
 Pick the cutoff function $\zeta(x)$ \! such that $\zeta \! = \! 1$  on  $K_{(1-\sigma)\frac{R_n}{2}},\sigma\in(0,1),\! |\nabla \zeta|\leq \frac{2}{\sigma R_n}.$
The  logarithmic  estimate  (Lemma \ref{lem:log-estimate})   with $ \bar Q= K_{\frac {R_n}2} \times (s,0]$ yields 
\begin{align}
 \label{logest}
&\underset{ s < t\leq 0}{\text{sup}} \int_{K_{(1-\sigma)\frac{R_n}{2}}}
\psi^2(x,t) \zeta^2 dx  \leq  \int_{K_{\frac{R_n}2}}  \psi^2(x,s) \zeta^2 dx \\ \notag
&+\gamma   \Big(\frac 1{c^2} \big[1+\log \frac H c\big] I_{\mathfrak d} (\mu_n^+)^{2{\mathfrak q}-m-1} |K_{\frac{R_n}{2}}|) \Big) 
+   \gamma  \log \frac Hc  \iint_{\bar Q}u^{m-1} |\nabla \zeta |^2dxdt \\
&=J_1+J_2+J_3.\notag
\end{align}
By \eqref{eq:psi-bounds} we have 
\begin{equation*}
\begin{aligned}
& J_1\,\leq  \big(\log \frac{H}{c}\big)^2|\{u(.,s) > \mu^+ -\frac{1}{4} \omega_n\}| \leq  \ (i \  \log 2)^2(1-\frac{1}{2}\nu)|K_{\frac{R_n}{2}}|.
\end{aligned}
\end{equation*}
From $ u-(\mu_n^+-\frac{\omega_n}{4})  >0$  and  \eqref{tr34},
we have
\begin{equation*}
\begin{aligned}
&u > \mu_n^+ - \frac{\omega_n}{4} >\frac 34 \omega_n - \frac{\omega_n}{4}= \frac{\omega_n}2.
\end{aligned}
\end{equation*}
Thus, taking in account the assumptions $ \mu_n^+\leq \frac{13}{12}\omega_n $ \,and\,$\omega_n\leq 1$, $ \theta_n= \omega_n^{1-m}$
we obtain
\begin{equation*}
\begin{aligned}
&J_2 \leq \gamma (\frac{13}{12}\omega_n)^{2{\mathfrak q} -m-1} \  \ \frac 1{c^2} [1+\log \frac H c \  I_{\mathfrak d} \  |K_{\frac{R_n}{2}}|
\\
&\leq \gamma_1  (\omega_n)^{2{\mathfrak q} -2m-2} [1+i \log 2] I_{\mathfrak d}  \
|K_{\frac{R_n}2 }|  
\leq\gamma_1 [1+i \log 2]I_{\mathfrak d}  |K_{\frac{R_n}{2}}|
\end{aligned}
\end{equation*} 
with $\gamma_1= \gamma (\frac{13}{12})^{2{\mathfrak q} -m-1 }$,
  and \,$\omega_n^{2{\mathfrak q} -2m-2}<1.$ \\
\begin{equation*}
\begin{aligned}
&J_3\leq \gamma  \ (i \ \log2)(\frac{3}{4}\omega_n)^{m-1}\frac{1}{(\sigma \frac{R_n}{2})^2} |K_{\frac{R_n}{2}}|\omega_n^{1-m}(\frac{R_n}{2})^2
\leq \gamma_2\frac{i}{\sigma^2}|K_{\frac{R_n}{2}}|.
\end{aligned}
\end{equation*} 
where \,$ \gamma_2 = \gamma \log 2 (\frac{3}{4})^{m-1}.$
Collecting these estimates  for $J_1,J_2,J_3$  we have
\begin{equation*}
\label{cr1}
\begin{aligned}
&\int_{K_{(1-\sigma)\frac{R_n}{2}} }\psi^2 (u)(x,t)dx\leq \,(i \  \log2)^2(1-\frac{1}{2}\nu)|K_{\frac {R_n}{2}}| \\
&+\gamma_1  [1+i \ \log 2 ]I_{\mathfrak d}  |K_{\frac{R_n}{2}}|+\gamma_2\frac{i}{\sigma^2}|K_{\frac{R_n}{2}}|.
\end{aligned}
\end{equation*}

We estimate below the left side  of \eqref{logest}  by integrating over a smaller set
$\hat S= \{u>\mu^+_n - \frac \omega{2^{i+2} }   \} \cap K_{(1-\sigma)\frac{R_n}{2}}$.  Observing that
$\psi$ is a decreasing function of $H$ in $\hat S$ and $H \leq \frac {\omega_n}4,$ we obtain\\
 $$ \psi \geq 
 \log^+ \Bigg[\frac{ \frac {\omega_n}4}
 { \frac {\omega_n}4
       -\Big[\mu^+_n-  \frac {\omega_n}{2^{i+2}} -(\mu_n^+ -\frac{\omega_n}{4})\Big]+c}\Bigg]= \log \ 2^{i-1}   , 
       $$
then
\begin{equation}
\label{cr2}
 \psi \geq (i-1)\ \log2.
\end{equation} 

Using \eqref{logest} and \eqref{cr2}  and for brevity setting  \,$ \hat k = \mu_n^+ -\frac{\omega_n}{2^{i+2}}$\,\, \\  one can find
\begin{align}
\label{cr3}
&((i-1) \ \log2)^2 |\{u(.,t)>\hat k\}\cap K_{(1-\sigma)\frac{R_n}{2}}|
\leq
(i \  \log2)^2(1-\frac{1}{2}\nu)|K_{\frac{R_n}{2}}| \\
&+ \gamma_1  [1+i\ \lg 2 ]I_{\mathfrak d}  |K_{\frac{R_n}{2}}|+\gamma_2\frac{i}{\sigma^2}| K_{\frac{R_n}{2}}|.\nonumber
\end{align}

 We remark that 
\begin{align}
\label{cr4}
&|\{u(\cdot,t)>\hat k\}\cap K_{\frac{R_n}{2}}|\! \leq |\{u(\cdot,t)>\hat k \}\cap K_{(1-\sigma)\frac{R_n}{2}})|+|K_{\frac{R_n}{2}}\! \setminus \! K_{(1-\sigma)\frac{R_n}{2}})| \\
&\notag \leq | \{u(\cdot,t)>\hat k\}\cap K_{(1-\sigma)\frac{R_n}{2}})|+|K_{\frac{R_n}{2}}|(1-(1-\sigma)^N ).\notag
\end{align}
Since \,$ (1-(1-\sigma)^N) <N\sigma$,\, from \eqref{cr4} we derive
\begin{equation}
\label{cr5}
\begin{aligned}
&-N\sigma|K_{\frac{R_n}{2}}|+|\{u(.,t)>\hat k\}\cap K_{\frac{R_n}{2}}|
<|\{u(.,t)>\hat k \}\cap K_{(1-\sigma)\frac{R_n}{2}}|.
\end{aligned}
\end{equation}
Inserting \eqref{cr5} in \eqref{cr3} gives
\begin{equation*}
\begin{aligned}
&((i-1) \ \log2)^2\big[|\{u(. ,t)>\hat k\}\}\cap K_{ \frac{R_n}{2}}|-N\sigma |K_{\frac{R_n}{2}}| \big]\\ 
&\, \leq (i\  \log2)^2(1-\frac{1}{2}\nu)|K_{\frac{R_n}{2}}| +\gamma_1  [1+i\    \lg( 2) ]I_{\mathfrak d}  | K_{\frac{R_n}{2}}|+
 \ \gamma_2\frac{i}{\sigma^2}|K_{\frac{R_n}{2}}|.
\end{aligned}
\end{equation*}
Dividing by \,$\,((i-1) \log 2)^2\,$ we obtain with  $\hat k= \mu_n^+ - \frac{ \omega_n} {2^{i+2} }   $
\begin{equation*}
\label{dis}
\begin{aligned}
&|\{u(. ,t)> \mu_n^+ - \frac{ \omega_n} {2^{i+2} } \}| \leq (\!1-\frac{\nu}{2})
 |K_{\frac{R_n}{2}}| \\
&
\times  \!\!\Bigg \{\! \Big[ \frac i{i-1} \Big]^2  + \frac{1}{( 1 -\frac{\nu}{2})}\Big[\gamma_1\frac{[1+i  \log 2 ] \ I_{\mathfrak d} }{((i-1) \log2)^2}+ \! \gamma_2 \frac{i}{\sigma^2}\frac{1}{((i-1)\  \log2)^2}+N\sigma\Big]\!
 \Bigg \}.
\end{aligned}
\end{equation*}
Choose  $ i $ large  and $\sigma $ small  in \eqref{dis} to obtain $$\frac{1}{(1-\frac{\nu}{2})} \Big[\gamma_1\frac{[1+i\  \lg2 ] \ I_{\mathfrak d} }{(i \ \log2)^2}+\gamma_2\frac{i}{\sigma^2}\frac{1}{(i \ log2)^2}+ N\sigma\Big] < \frac{\nu}{2},$$
so that  \eqref{nu*} holds with 
\begin{equation}
\label{n*}
 n_* = i+2.
\end{equation} 
 \end{proof}
Now let $\nu, n_* $ be the numbers determined in  Lemma \ref{n*nu} and \eqref {n*} respectively. \\
As a consequence we have that
\begin{equation}
\label{B.10.4}
|\{\{u(,t) > \mu_n^+ - \frac{\omega_n}{2^j}\}\cap K_{\frac{R_n}{2}}| <(1-\frac{1}{4}\nu^2)|K_{\frac{R_n}{2}}|
\end{equation}
 for all  $j \geq n_*$ 
and for all times
\begin{equation}
\label{B.10.5}
 -\frac{1}{2}\nu \omega_n^{1-m}\Big(\frac{R_n}{2}\Big)^2 <t< 0.
\end{equation}
Now let us introduce the cylinder 
\begin{equation*}
\label{B.11.1}
Q_{\frac{R_n}{2}}(\theta_*) = K_{\frac{R_n}{2}}\times(-\theta_*(\frac{R_n}{2})^2,0],
\end{equation*}
 with 
 
 \begin{align}\label{theta*}
  \theta_* =\frac{1}{2}\nu\omega_n^{1-m}.
  \end{align}

From the informations obtained in \eqref{B.10.4}  -\eqref{B.10.5} we can deduce that in the cylinder $ Q_{\frac{R_n}{2}}(\theta_*) $
 the set where $u$ is close to its supremum $\mu_n^+$ can be made sufficiently small. Firstly  we establish an estimate of $ \iint_{Q_{\frac{R_n}{2}}(\theta_*)}|\nabla(u-k_j) |^2 dxdt. $


\begin{lemma}
\label{grad u}
Let
  $n_*$ defined in \eqref{n*}, $j\geq n_*$.
There exists a positive constant $\bar \gamma$ such that
\begin{equation}
\label{B.11.3}
\iint_{Q_{\frac{R_n}{2}}(\theta_*)}|\nabla(u-k_j) |^2 dx dt\leq \frac{\bar{\gamma}}{\nu( \frac{R_n}{2})^2}\big(\frac{\omega_n}{2^j}\big)^2 |Q_{\frac{R_n}{2}}(\theta_*)|,
\end{equation}
with $\theta_{*}$ defined in \eqref{theta*}.
 \end{lemma}
 \begin{proof}
  Taking into account \eqref{B.10.5},  construct the new cylinder
   $ Q'= K_{R_n}\times(-\nu\omega_n^{1-m}(\frac{R_n}{2})^2, \ 0]$. Clearly $Q_{\frac{R_n}{2}}(\theta_*) \subset Q'$.
We recall that  $\frac 34 \omega_n < \mu_n^+< \frac{13 } { 12}\omega_n $. 
Let us consider the truncated functions 
\begin{equation}
 \label{q*}
\begin{aligned}
 u-k_j  \quad   \text{with} \quad  k_j = \mu_n^+ -\frac{\omega_n}{2^j},
 \quad j= n_*,... , n_* + q_* 
\end{aligned}
\end{equation}  
$q_* $ to be determined
and  apply   the energy estimates in  Lemma \ref{lem:energy2}.
where we  neglect the first term on the left.
 The cutoff function $\eta=1$   on $ Q_{\frac{R_n}{2}}(\theta_*) $ and $\eta=0$ on the parabolic boundary of  $ Q' $ with
 $$ |\nabla \eta|< \frac{2}{\frac {R_n}2},\qquad  0\leq \eta_t \leq \frac{8}{\nu\theta R_n^2}\quad \theta = \omega^{1-m}.$$
 With these assumptions \eqref{EEu>k} takes the form  
\begin{align}
 \label{B.11.a}
&\iint_{Q'}u^{m-1}|\nabla  ((u-k_j) \eta)|^2 dx d\tau \\
& \notag \leq \gamma \Big(\iint_{ Q'} u^{m-1}(u-k_j)^2|\nabla \eta|^2 dxd\tau + I_{\mathfrak d}  \ \omega_n^{2{\mathfrak{q}-m-1}}|Q'| \\   \nonumber
&+ \int_{K_{R_n}}(u-k_j)^2\eta^2(x,- \nu \omega_n^{1-m} ( \frac{R_n}{2})^2) +\iint_{Q'}  ( u-k_j)^2      \eta  \eta_t  \ dx  dt\Big).\\  \nonumber
\end{align}

By the choice of $\eta, $ the third term on the right hand of the last inequality vanishes.
Since $ u>k_j $ (with $ k_j =\mu_n^+ -\frac{\omega_n}{2^j}$ and  \,$ \mu_n^+> \frac{3}{4}\omega_n$)  we have

$$ u > \mu_n^+ -\frac{\omega_n}{2^j}>\frac{3}{4}\omega_n-\frac{\omega_n}{2^j}\geq \frac{1}{4}\omega_n.$$
Moreover $ u \leq \mu_n^+< \frac{13}{12}\omega_n.$ \ \ 
From \eqref{B.11.a}  we derive, after dividing by $\omega_n^{m-1}$
 \begin{align}
 \label{B.11.b}
&\big(\frac{13}{12} \big)^{m-1}\! \iint_{Q}|\nabla(u-k_j) |^2 dx d\tau \\ 
& \notag \leq \gamma \Bigg(\big(\frac{1}{4}\big)^{m-1}\frac{1}{(\frac{R_n}{2})^2} \big(\frac{\omega_n}{2^j}\big)^2 |Q'| +
 \underset{drift \ term }   {\Big\{ I_{\mathfrak d}  \ \omega_n^{2{\mathfrak{q}-2m}}|Q'|\Big \}}+\!
\big(\frac{\omega_n}{2^j}\big)^2  \! \frac{1}{\nu (\frac{R_n}{2})^2}|Q'|  \Bigg). \\  \nonumber
\end{align}

%
{\bf Estimate of the drift term}\\
Let us reduce the drift term
$
\Big\{I_{\mathfrak d}  \ \omega_n^{2{\mathfrak{q}-2m}}|Q_{\frac{R_n}{2}}(\theta_*)| \Big\}
$
to a term similar to \\
$I_{\mathfrak d}   \big(\frac{\omega}{2^j}\big)^2  \ \frac{1}{\nu (\frac{R_n}{2})^2}|Q_{\frac{R_n}{2}}(\theta_*)|. $
 Multipling and dividing by $\nu(\frac{R_n}{2})^2 
(2^j)^2$ , we have
\begin{equation*}
I_{\mathfrak d}  \omega_n^{2{\mathfrak{q}-2m}}|Q_{\frac{R_n}{2}}(\theta_*)|  \leq I_{\mathfrak d}  \ \omega_n^{2{\mathfrak{q}-2m}}   \ \frac{ \omega_n^2  }{  \nu (\frac{R_n}{2})^2 \ (2^j)^2}  |Q_{\frac{R_n}{2}}(\theta_*)|
\end{equation*} 
where $\nu<1$,  from \eqref{q*}
 \ $j=n_*,..., n_*+q_*$ and by choosing in  \eqref{lambda}  $\lambda= n_*+q_*$, then 
 $\frac{R_n}{2} 2^j < \frac{R_n}{2} \ 2^{n_*+q_*} <\omega;$    moreover  $\omega_n^{2{\mathfrak{q}-2m}} <1.$  
 For a suitable constant $ \bar{\gamma}$  \eqref {B.11.b} becomes
 \eqref{B.11.3}.\end{proof}
Now we are in the position to derive an estimate of the set where $u$ is near to $\mu_n^+$ within $|Q_{\frac{R_n}{2}}(\theta_*)|$.


\begin{lemma}
\label{grad u2}
 In the hypotheses of Lemma \ref{grad u} , for every  $\bar  \nu_*\in (0,1)$ there exists $ q_*,$ depending only on the data and $\bar  \nu_*, $ such that 
 \begin{equation*}
 \label{B.11.2}
 |\{u > \mu_n^+ - \frac{\omega}{2^{n_* + q_*}}\}\cap Q_{\frac{R_n}{2}}(\theta_*)|\leq \bar \nu_*|Q_{\frac{R_n}{2}}(\theta_*)|.
 \end{equation*}
 \end{lemma}
 
 \begin{proof}
 We observe that
 as a consequence of  \eqref{B.10.4}  we have 
\begin{equation*}
|\{u(.,t)< k_j \}\cap K_{\frac{R_n}{2}}|\geq\frac{1}{4} \nu^2|K_{\frac{R_n}{2}}|.
\end{equation*} 

We start with an application of the DeGiorgi theoretical Lemma \ref{lemma:DeGiorgi} of Section \ref{Pre} to the function $u(x,t)$ in the range \eqref{B.10.5}
over the cube $K_{\frac{R_n}{2}}$ for the levels
\begin{equation*}
k = k_j < \ell= k_{j+1} \quad \text{so\, that} \quad (\ell-k) = (k_{j+1} -k_j) =\frac{\omega}{2^{j+1}},\ \ \ 
 k_j =\mu_n^+ -\frac{\omega_n}{2^j}.
\end{equation*}
Consequently
\begin{align}
\label{B.11.3bis}
&\frac{\omega_n}{2^{j+1}} \big|\{u(\cdot,t)> k_{j+1}\}\cap K_{\frac{R_n}2 }\big|  
\\ \nonumber
&\leq
 \frac { \gamma_{D} (\frac{R_n}2)^{N+1}   } 
  {| \{u(.,t)< k_{j} \}  \cap K_{\frac{R_n}{2}}|} \int_{\{k_ j< u< k_{j+1}\}\cap K_{\frac{R_n}{2} } } |\nabla u| dx \\  \nonumber
& \leq 4 \gamma_{D}  \frac
{\frac{R_n}{2} }  {\nu^2}\!
\! \Big( \int_{ \{k_ j< u< k_{j+1}\}\cap K_{\frac{R_n}{2} }  } \!\! \!
|\nabla u(\cdot, t)|^2 dx\Big)^{\frac{1}{2}} \!  \\  \nonumber 
&  \times
 \Big|(\{u(.,t) >k_j\} -\{u(.,t)> k_{j+1}\! \})\! \cap \! K_{ \frac{R_n}{2}} \Big|^{\frac{1}{2}}, \\  \nonumber 
\end{align}

where  in  \eqref{B.11.3bis} the H\"{o}lder inequality is used.\\
Set
 $A_j =\{u > k_j\}\cap Q_{\frac{R_n}{2}}(\theta_*)  \  {\text {with}}  \  |A_j|= \int_{-\theta_* (\frac{R_n}{2})^2}^{0} \ |\{u(\cdot,t)> k_j\}\cap K_{\frac{R_n}{2}}| dt.$\\
 Integrate  the inequality \eqref{B.11.3bis} over the range $(- \theta_*  ( \frac{R_n}{2} )^2, 0]$  to obtain
\begin{equation*}
\label{B.11.3tris}
\begin{aligned}
&\frac{\omega_n}{2^{j+1}}|A_{j+1}|  \leq 4 \gamma_{D} \frac{\frac{R_n}{2}}{\nu^2}\Big(\iint_{Q_{\frac{R_n}{2}}(\theta_*)} |\nabla (u-k_j)|^2dx dt\Big)^{\frac{1}{2}}
\big(| A_j|- |A_{j+1}|\big)^{\frac{1}{2}}  .
\end{aligned}
\end{equation*}
Square both sides  and apply  the estimate \eqref{B.11.3} in Lemma \ref{grad u}
 for the term  containing $|\nabla (u-k_j)|$   to obtain with $\gamma_2 = \bar \gamma (4 \gamma_{D})^2 $
$$|A_{j+1}|^2 \leq \  \frac{\gamma_2}{\nu^5} \ |Q_{\frac{R_n}{2}}(\theta_*)|\big( |A_j| -|A_{j+1}|\big).$$
For $ j = n_* +1, n_*+2,..., (n_* + q_* -1),$ we
sum  the previous  inequalities and  
since $ |A_{j+1} |\geq |A_{j+2}| \geq...\geq |A_{n_* +q_*}|,$ we have\\ 
$  (q_* -2) |A_{n_* +q_* }|^2  \leq \frac{\gamma_2}{\nu^5}|Q_{\frac{R_n}{2}}(\theta_*)| \sum\limits_{j=n_* +1} ^{n_* + q_* -1}  |A_{j}|- |A_{j+1}| \leq \frac{\gamma_2}{\nu^5}
|Q_{\frac{R_n}{2}}(\theta_*)|^2. $ \\
It follows
\begin{equation}
\label{B.11.4}
|A_{n_* +q_* }|\leq  \sqrt{\frac{\gamma_2}{({q_* -2})\nu^5}} \ |Q_{\frac{R_n}{2}}(\theta_*)|=: 
 \bar \nu_*|Q_{\frac{R_n}{2}}(\theta_*)|.
\end{equation} 
Once fixed $\bar \nu_* $  we compute $ q_*$ from  \eqref{B.11.4} and  Lemma \ref{grad u2} is proved.  
\end{proof}


Subsequent arguments (energy estimates, measure decay of superlevel sets, and a De Giorgi-type lemma) yield the oscillation decay
\begin{equation}
\label{2osc}
u \leq \mu_n^+ - \tfrac{\omega_n}{2^{n_*+q_*+1}} \quad \text{in } Q_{\frac {R_n}4}(\theta_*), 
\end{equation}
which shows that $u$ is strictly below its supremum $\mu^+$  in a smaller cylinder.\\

 \subsection{Proof of  Proposition \ref{prop}  and  Theorem~\ref{thm:holder} }

We now turn to the proof of H\"older continuity by combining the two alternatives.\\
  We begin by assuming that the condition \eqref{umagg} of the first alternative does not hold. In the case of the second alternative, we already obtained inequality \eqref{2osc}, which holds in the cylinder $Q_{\frac {R_n}4}(\theta_*)$. Consequently, under hypothesis \eqref{tr34}, we deduce
\begin{equation}\label{sup u}
\operatorname*{ess\,sup}_{\frac {R_n}4 (\theta_*)} u \leq \mu^{+}_n - \frac{\omega_n}{2^{q_*+n_*+1}}.
\end{equation}

As already remarked in inequality (8.7), the left-hand side can always be assumed, while the right-hand side holds provided that 
\[
\mu^{-}_n < \frac{\omega}{12},
\] 
which is consistent with (6.1), namely $\omega \geq \tfrac{12}{13} \mu^{+}$ \  and thus the second alternative is valid. Subtracting  \ 
$
\operatorname*{ess\,inf}_{Q_{\frac {R_n}4}(\theta_*)} u > \mu^{-}_n
$
\ from both sides of \eqref{sup u}, we obtain
\[
\operatorname*{ess\,osc}_{Q_{\frac {R_n}4}(\theta_*)} u \leq \Big(1 - \frac{1}{2^{q_*+n_*}}\Big)\,\omega_n.
\]

By construction, the cylinder $Q_{n+1}(\theta)$ is given by
\[
Q_{n+1}(\theta) = K_{R_{n+1}} \times \Big(-\omega_{n+1}^{1-m} R_{n+1}^2, 0 \Big],
\quad R_{n+1} = \frac{R_n}{b}, \quad b > 1,
\]
while
\[
Q_{\frac {R_n}4}(\theta^{\ast}) = K_{\frac {R_n}4} \times \Big(-\tfrac{1}{2}\nu \omega_n^{1-m} ({\frac {R_n}4})^2, 0 \Big],
\quad \theta_* = \tfrac{1}{2}\nu \omega_n^{1-m}.
\]
Choosing $b = \sqrt{32/\nu}$, we guarantee the inclusion
$
Q_{n+1}(\theta) \subset Q_{\frac {R_n}4}(\theta_*).
$
Thus,
\begin{equation*}\label{eq:833}
\omega_{n+1} := \operatorname*{ess\,osc}_{Q_{n+1}(\theta)} u 
\leq \operatorname*{ess\,osc}_{Q_{\frac {R_n}4}(\theta_*)} u 
\leq \delta \omega_n,
\qquad \delta = 1 - \frac{1}{2^{q_*+n_*+1}}.
\end{equation*}
Hence, in the second alternative, the oscillation of $u$ decreases when passing to a smaller cylinder.\\ 
We now show that the same conclusion holds in the first alternative. In this case,\eqref{umaggg}
 and \eqref{tr34} imply
$$
-\operatorname*{ess\,inf}_{Q_{\frac {R_n}4}(\theta)} u < -\mu^{-}_n - \frac{\omega_n}{4}.
$$
Adding 
$\underset{Q_{\frac {R_n}4}(\theta)}{\operatorname*{ess\,sup}\ u}$  
 yields
\[
\operatorname*{ess\,osc}_{Q_{\frac {R_n}4}(\theta)} u 
\leq \mu^{+}_n - \mu^{-}_n - \frac{\omega_n}{4}
= \frac{3}{4}\omega_n.
\]
Since, by construction, $Q_{n+1}(\theta) \subset Q_{\frac {R_n}4}(\theta)$ (recall that $1/b = \sqrt{\nu/32} < 1/4$), we deduce
\[
\omega_{n+1} = \operatorname*{ess\,osc}_{Q_{n+1}(\theta)} u
\leq \operatorname*{ess\,osc}_{Q_{\frac {R_n}4}(\theta)} u
\leq \frac{3}{4}\omega_n.
\]

Finally, we remark that $\tfrac{3}{4} < 1 - 2^{-(q_*+n_*+1)}$, so that in both alternatives we obtain the uniform decay
\[
\omega_{n+1} \leq \delta \,\omega_n, \qquad \delta \in (0,1).
\]

This completes the ``Induction Argument'' introduced in Proposition \ref{prop}. As a consequence, the oscillation of $u$ decreases geometrically, and therefore $u$ is H\"older continuous. As noted in the introduction, the regularity of $v$ follows in a straightforward manner (see, e.g., \cite{L},\cite{MRVVp}).

\begin{remark}[Perspectives]
The above argument establishes H\"older regularity in the range $ \tfrac{(N-2)_{+}}{N+2}<m<1$. In future work, we shall investigate the critical and subcritical regimes, namely $0 < m \leq \tfrac{(N-2)_{+}}{N+2}$, as well as the borderline logarithmic case. These settings are expected to require new ideas and refined techniques, since the oscillation decay mechanism becomes more delicate near the critical exponent.
\end{remark}

\section*{Acknowledgments}

This work was carried out within the activities of the 
\textit{Gruppo Nazionale per l'Analisi Matematica, la Probabilit$\grave a$ e le loro Applicazioni} (GNAMPA) 
of the \textit{Istituto Nazionale di Alta Matematica} (INdAM). 
The authors M.M., S.V-P., and V.V. gratefully acknowledge this support.\\
M. Marras is partially supported by the research project "Partial Differential Equations and their role in understanding natural phenomena", \\
CUP F23C25000080007, funded by Fondazione Banco di Sardegna annuality (2023); by the project National Recovery and Resilience Plan (NRRP), Mission 4 Component 2 Investment 1.5 - Call for tender No.3277 published on December 30, 2021 by the Italian Ministry of University and Research (MUR) funded by the European Union-NextGenerationEU. Project Code ECS0000038-Project Title eINS Ecosystem of Innovation for Next Generation Sardinia-CUP F53C22000430001- Grant Assignment Decree No. 1056 adopted on June 23, 2022 by the Italian Ministry of University and Research (MUR) and and by the grant   INDAM-GNAMPA Project, CUP E53C25002010001. \\

\end{document}